\documentstyle[12pt,a4,amsfonts]{article}

\newtheorem{Theorem}{Theorem}[section]

\newtheorem{Lemma}[Theorem]{Lemma}

\newtheorem{Remark}[Theorem]{Remark}
\newtheorem{Example}[Theorem]{Example}

\newtheorem{Claim}[Theorem]{Claim}
\newtheorem{Problem}[Theorem]{Problem}

\def\RMN#1{\uppercase\expandafter{\romannumeral#1}}

\newcommand{\specialarrow}[3]{\mathrel{\mathop{#1}\limits^{#2}_{#3}}}

\newcommand{\uetuki}[2]{{\mathop{\scriptstyle {#1}}\limits_
{\scriptstyle {#2}}}}

\newcommand{\qed}{{\unskip\nobreak\hfil\penalty50\quad\null\nobreak\hfil{\bf
q.e.d.}\parfillskip0pt\finalhyphendemerits0\par\medskip}}
\newcommand{\proof}{\noindent{\it Proof.} \ }

\newcommand{\rank}{\mathop{\rm rank}\nolimits}

\newcommand{\spec}{\mathop{\rm Spec}\nolimits}
\newcommand{\proj}{\mathop{\rm Proj}\nolimits}

\newcommand{\subq}{_{\Bbb Q}}

\begin{document}

\title{
A local ring such that 
the map \\ between Grothendieck groups with rational coefficient \\
induced by completion is not injective}
\author{Kazuhiko Kurano and Vasudevan Srinivas}
\date{}
\maketitle

\abstract{In this paper, we construct a local ring $A$ such that
the kernel of the map $G_0(A)\subq \rightarrow G_0(\hat{A})\subq$ 
is not zero, where $\hat{A}$ is the comletion of $A$ with respect to 
the maximal ideal, and $G_0( \ )\subq$ is the Grothendieck group 
of finitely generated modules with rational coefficient.
In our example, $A$ is a two-dimensional local ring which is essentially of finite type over
${\Bbb C}$, but it is not normal.}

\section{Introduction}
For a Noetherian ring $R$, we set
\[
G_0(R) = 
\frac{
\bigoplus\limits_{\mbox{\scriptsize $M$: f.\ g.\ $R$-mod.\ }} {\Bbb Z} [M] 
}{
\left< [L] + [N] - [M] \left|
\mbox{$ 0 \rightarrow L \rightarrow M \rightarrow N \rightarrow 0 $ is exact}
\right. \right> ,
}
\]
that is called the {\em Grothendieck group} 
of finitely generated $R$-modules.
Here, $[M]$ denotes the free basis (corresponding to 
a finitely generated $R$-module $M$) of the free module 
$\bigoplus {\Bbb Z}[M]$,
where ${\Bbb Z}$ is the ring of integers.

For a flat ring homomorphism $R \rightarrow A$,
we have the induced map $G_0(R) \rightarrow G_0(A)$ 
defined by $[M] \mapsto [M \otimes_RA]$.

We are interested in the following problem (Question~1.4 in \cite{KK}):
\begin{Problem}\label{prob1}
\begin{rm}
Let $R$ be a Noetherian local ring.
Is the map 
$G_0(R)_{\Bbb Q} \rightarrow G_0(\widehat{R})_{\Bbb Q}$ injective?
\end{rm}
\end{Problem}
Here, $\widehat{R}$ denotes the ${\frak m}$-adic completion of $R$,
where ${\frak m}$ is the unique maximal ideal of $R$.
For an abelian group $N$, $N\subq$ denotes the tensor product with 
the field of rational numbers ${\Bbb Q}$.

We shall explain motivation and applications.

Assume that $S$ is a regular scheme and $X$ is 
a scheme of finite type over $S$.
Then, by the singular Riemann-Roch theorem~\cite{F}, 
we obtain an isomorphism 
\[
\tau_{X/S} : G_0(X)_{\Bbb Q} \stackrel{\sim}{\rightarrow} 
A_*(X)_{\Bbb Q},
\]
where $G_0(X)$ (resp.\ $A_*(X)$) is the {\em Grothendieck group}
of coherent sheaves on $X$ (resp.\ {\em Chow group} of $X$).
We refer the reader to Chapters~1, 18, 20 in \cite{F} for definition 
of $G_0(X)$, $A_*(X)$ and $\tau_{X/S}$.
Note that $G_0(X)$ (resp.\ $\tau_{X/S}$) 
is denoted by $K_0(X)$ (resp.\ $\tau_{X}$) in \cite{F}.
The map $\tau_{X/S}$ usually depends on the choice of $S$.
In fact, we have
\begin{eqnarray*}
\tau_{{\Bbb P}^1_k/{\Bbb P}^1_k}([{\cal O}_{{\Bbb P}^1_k}]) & = & 
[{\Bbb P}^1_k] 
\in A_*({\Bbb P}^1_k)_{\Bbb Q} = {\Bbb Q} [{\Bbb P}^1_k] \oplus {\Bbb Q} [t]
\\
\tau_{{\Bbb P}^1_k/\spec k}([{\cal O}_{{\Bbb P}^1_k}]) & = & 
[{\Bbb P}^1_k] + \chi({\cal O}_{{\Bbb P}^1_k}) [t] =
[{\Bbb P}^1_k] +  [t]
\in A_*({\Bbb P}^1_k)_{\Bbb Q} ,
\end{eqnarray*}
where $t$ is a $k$-rational closed point of ${\Bbb P}^1_k$ over a field $k$.
Here, for a closed subvariety $Y$, $[Y]$ denotes the algebraic cycle
corresponding to $Y$.
Hence,
\[
\tau_{{\Bbb P}^1_k/{\Bbb P}^1_k}([{\cal O}_{{\Bbb P}^1_k}]) \neq
\tau_{{\Bbb P}^1_k/\spec k}([{\cal O}_{{\Bbb P}^1_k}])
\]
in this case.
However, for a local ring $R$ which is a homomorphic image of a regular
local ring $T$, the map $\tau_{\spec R/\spec T}$ is 
independent of the choice of $T$ in many cases.
In fact, if $R$ is a complete local ring
or $R$ is essentially of finite type over either a field or 
the ring of integers, it is proved in Propopsition~1.2 of \cite{K16} that
the map $\tau_{\spec R/\spec T}$ 
is actually independent of $T$.

From now on, for simplicity, we denote $\tau_{\spec R /\spec T }$
by $\tau_{R/T}$.
It is natural to ask the following:

\begin{Problem}\label{prob2}
\begin{rm}
Let $R$ be a homomorphic image of a regular local ring $T$.
Is the map $\tau_{R/T}$ independent of $T$?
\end{rm}
\end{Problem}

Remark that, by the singular Riemann-Roch theorem, 
the diagram
\[
\begin{array}{ccc}
G_0(R)_{\Bbb Q} & \stackrel{\tau_{R/T}}{\longrightarrow} & 
A_*(R)_{\Bbb Q} \\
\downarrow & & \downarrow \\
G_0(\widehat{R})_{\Bbb Q} & \stackrel{\tau_{\widehat{R}/\widehat{T}}}{\longrightarrow} & 
A_*(\widehat{R})_{\Bbb Q}
\end{array}
\]
is commutative,
where the vertical maps are induced by the completion
$R \rightarrow \widehat{R}$.
We want to emphasize that the bottom map, as well as the vertical maps is
independent of the choice of $T$
since $\widehat{R}$ is complete (Propopsition~1.2 of \cite{K16}).
Therefore, if the vertical maps are injective,
then the top map is also independent of $T$.

Therefore, if Problem~\ref{prob1} is affirmative, then so is Problem~\ref{prob2}.

\vspace{3mm}

We shall explain another motivation.

Roberts~\cite{R2} and Gillet-Soul\'e~\cite{GS} 
proved the vanishing theorem of intersection 
multiplicities for complete intersections.
If a local ring $R$ is a complete intersection,
then $\tau_{R/T}([R]) = [\spec R ]$ holds, where
\[
[\spec R ] = \sum_{\uetuki{{\frak p} \in \spec R }{\dim R/{\frak p} = \dim R}} 
\ell_{R_{\frak p}}(R_{\frak p}) [\spec R/{\frak p} ]
\in A_{\dim R}(R)\subq .
\]
In \cite{R2}, Roberts proved the vanishing theorem of 
intersection multiplicities
not only for complete intersections
but also for local rings satisfying $\tau_{R/T}([R]) = [\spec R ]$.
Inspired by his work, 
Kurano~\cite{K16} started to study local rings
which satisfy the condition $\tau_{R/T}([R]) = [\spec R ]$,
and call them {\em Roberts} rings.
If $R$ is a Roberts ring, then the completion, the henselization and
localizations of it are also Roberts rings~\cite{K16}.
However, the following problem remained open.

\begin{Problem}\label{prob3}
\begin{rm}
If $\widehat{R}$ is a Roberts ring, is $R$ so?
\end{rm}
\end{Problem}

It is proved in Proposition~6.2 of \cite{KK} that 
Problem~\ref{prob3} is affirmative
if and only if so is Problem~\ref{prob1}.

\vspace{2mm}

The following partial result on Problem~\ref{prob1} was given by
Theorem~1.5 in \cite{KK}:

\begin{Theorem}[Kamoi-Kurano, 2001~\cite{KK}]\label{KK}
Let $R$ be a homomorphic image of an excellent regular local ring.
Assume that $R$ satisfies 
one of the following three conditions:

(i) $R$ is henselian,

(ii) $R = S_{\frak n}$, 
where $S$ is a standard graded ring over a field 
and ${\frak n} = \oplus_{n > 0}S_n$,

(iii) $R$ has only isolated singularity.

Then, the induced map
$G_0(R) \rightarrow G_0(\widehat{R})$ is injective.
\end{Theorem}

However, the following example was given by Hochster:

\begin{Example}[Hochster \cite{H}]\label{hochster}
Let $k$ be a field. 
We set 
\begin{eqnarray*}
T & = & k[x,y,u,v]_{(x,y,u,v)} , \\
P & = & (x,y) , \\
f & = & xy - ux^2 - vy^2 .
\end{eqnarray*}

Then, 
${\rm Ker}(G_0(T/fT) \rightarrow G_0(\widehat{T/fT})) \ni [T/P] \neq 0$.
In this case, $2 \cdot [T/P] = 0$.
\end{Example}

The ring $T/fT$ is not normal in the above example.
Recently Dao~\cite{Dao} found the following example:

\begin{Example}[Dao \cite{Dao}]\label{dao} 
We set 
\begin{eqnarray*}
R & = & {\Bbb Q}[x,y,z,w]_{(x,y,z,w)}/(x^2 + y^2 - (w+1)z^2) , \\
P & = & (x,y,z) .
\end{eqnarray*}

Then, 
${\rm Ker}(G_0(R) \rightarrow G_0(\widehat{R})) \ni [R/P] \neq 0$.
In this case, $2 \cdot [R/P] = 0$.
Here, $R$ is a normal local ring.
\end{Example}

The following is the main theorem of this paper:

\begin{Theorem}\label{maintheorem}
There exists a $2$-dimensional local ring $A$,
which is essentially of finite type over ${\Bbb C}$, that satisfies
\[
{\rm Ker}(G_0(A)_{\Bbb Q} \rightarrow G_0(\widehat{A})_{\Bbb Q}) 
\neq 0 .
\]
\end{Theorem}

\begin{Remark}
\begin{rm}
\begin{enumerate} 
\item
By Theorem~\ref{maintheorem}, we know that both
Problem~\ref{prob1} and Problem~\ref{prob3} are negative.
That is to say, there exists a local ring $R$ such that 
$\widehat{R}$ is a Roberts ring,
but $R$ is not so.
\item
Problem~\ref{prob2} is still open.
\item
In \cite{K23}, we defined notion of numerical equivalence
on $G_0(R)$ and $A_*(R)$.
We set $\overline{G_0(R)} = G_0(R)/ \sim_{\rm num.}$
and $\overline{A_*(R)} = A_*(R)/ \sim_{\rm num.}$.
Then, we have the following:
\begin{enumerate}
\item
$\overline{G_0(R)} \rightarrow \overline{G_0(\widehat{R})}$ is injective
for any local ring $R$.
\item
The induced map $\overline{\tau_{R/T}} : \overline{G_0(R)}_{\Bbb Q} \stackrel{\sim}{\rightarrow} \overline{A_*(R)}_{\Bbb Q}$ 
is independent of $T$.
\item
$R$ is a numerically Roberts ring iff so is $\widehat{R}$.
(Definition of numerically Roberts rings was given in \cite{K23}.
The vanishig theorem of intersection multiplicities holds
true for numerically Roberts rings.)
\end{enumerate}
\item
The ring $A$ constructed in the main theorem is not normal.
We do not know any example of a normal local ring 
that does not satisfy Problem~\ref{prob1}.
\end{enumerate}
\end{rm}
\end{Remark}

Theorem~\ref{maintheorem} immediately follows from the following
two lemmas.

\begin{Lemma}\label{lem1}
Let $K$ be an algebraically closed field, and let
$S = \oplus_{n \geq 0}S_n$ be a standard graded ring over $K$,
that is, a Noetherian graded ring generated by $S_1$ over 
$S_0 = K$.
We set $X = \proj S $, and assume that $X$ is smooth over $K$
with $d = \dim X \geq 1$.
Let $h$ be the very ample divisor on $X$ of this embedding.
Let $\pi : Y \rightarrow \spec S $ be the blow-up at 
${\frak n} = \oplus_{n > 0}S_n$.

Assume the following:
\begin{enumerate}
\item
Set $R = S_{{\frak n}}$ and let $\widehat{R}$ be the completion of $R$.
Then, the map 
$A_1(R)_{\Bbb Q} \rightarrow
A_1( \widehat{R} )_{\Bbb Q}$ induced by completion is an isomorphism.
\item
There exists a smooth connected curve $C$ in $Y$ that satisfies
following two conditions:
\begin{itemize}
\item[(i)]
$C$ transversally intersects with $\pi^{-1}({\frak n}) \simeq X$ 
at two points, namely $P_1$ and $P_2$.
\item[(ii)]
$[P_1] - [P_2] \neq 0$
in $A_0(X)_{\Bbb Q}/h \cdot A_1(X)_{\Bbb Q}$.
\end{itemize}
\end{enumerate}

Then, there exists a local ring $A$ of dimension $d + 1$, which is
essentially of finite type over $K$, such that
\[
{\rm Ker}(G_0(A)_{\Bbb Q} \rightarrow G_0(\widehat{A})_{\Bbb Q}) 
\neq 0 .
\]
\end{Lemma}

\begin{Lemma}\label{lem2}
We set $S = {\Bbb C}[x_0, x_1, x_2]/(f)$, 
where $f$ is a homogeneous cubic polynomial.
Assume that $X = \proj S $ is smooth over ${\Bbb C}$.

Then, $R$ satisfies the assumption in Lemma~\ref{lem1} with $d = 1$.
\end{Lemma}

We shall prove the above two lemmas in the following sections.

\section{A proof of Lemma~\ref{lem1}}
Here, we shall give a proof of Lemma~\ref{lem1}.

Let ${\frak p}$ be the prime ideal of $S$ that satisfies 
$\spec S/{\frak p}  = \pi(C)$.
Set $R = S_{{\frak n}}$ and ${\frak m} = {\frak n}R$.

Then, $C$ is the normalization of $\spec S/{\frak p} $.
We denote by $v_i$
the normalized valuation of the discrete valuation ring at $P_i \in C$ 
for $i = 1, 2$.

First of all, we shall prove the following:

\begin{Claim}\label{claim1}
There exists $s \in {\frak m}/{\frak p}R$ such that
\begin{enumerate}
\item
$v_1(s) = v_2(s) > 0$, and
\item
$K[s]_{(s)} \hookrightarrow R/{\frak p}R$ is finite.
\end{enumerate}
\end{Claim}

\proof
Let $C'$ be the smooth projective connected curve over $K$ 
that contains $C$ as a Zariski open set.
We regard $P_1$, $P_2$ as points of $C'$.

Let $R(C')$ be the field of rational functions 
on $C'$.
Since $P_1$ is an ample divisor on $C'$,
there exists $t_1 \in R(C')^\times$ such that
\begin{itemize}
\item
$P_1$ is the only pole of $t_1$, and
\item
$P_2$ is neither a zero nor a pole of $t_1$.
\end{itemize}
Similarly, one can find $t_2 \in R(C')^\times$ such that
\begin{itemize}
\item
$P_2$ is the only pole of $t_2$, and
\item
$P_1$ is neither a zero nor a pole of $t_2$.
\end{itemize}
Replacing $t_1$ (resp.\ $t_2$) with a suitable powers of $t_1$ (resp.\ $t_2$),
we may assume $v_1(t_1) = v_2(t_2) < 0$.

Put $t = 1/t_1t_2 \in R(C')^\times$.
Then, $\{ P_1, P_2 \}$ is the set of zeros of $t$.
Note that $v_1(t) = v_2(t) > 0$.

Let $O_{v_i}$ be the discrete valuation ring at $P_i$ 
for $i = 1, 2$.
Then, $K[t]_{(t)}$ is a subring of 
\[
O_{v_1} \cap O_{v_2} = \overline{S/{\frak p}} \otimes_{S/{\frak p}} R/{\frak p}R ,
\]
where $\overline{( \ \ )}$ is the normalization of the given ring.

Since $\{ P_1, P_2 \}$ is just the set of zeros of $t$,
$O_{v_1} \cap O_{v_2}$ is the integral closure
of $K[t]_{(t)}$ in $R(C')$.
In particular, $\overline{S/{\frak p}} \otimes_{S/{\frak p}} R/{\frak p}R$ is finite
over $K[t]_{(t)}$.

Let $I$ be the conductor ideal of the normalization
\[
R/{\frak p}R \subset \overline{S/{\frak p}} \otimes_{S/{\frak p}} R/{\frak p}R .
\]
Let ${\frak m}_i$ be the maximal ideal of 
$\overline{S/{\frak p}} \otimes_{S/{\frak p}} R/{\frak p}R$
corresponding to $P_i$
for $i = 1, 2$.
Since $I$ is contained in ${\frak m}/{\frak p}R$, 
\[
I \subset {\frak m}_1 \cap {\frak m}_2 .
\]
Therefore, we have
\[
\sqrt{I} = {\frak m}_1 \cap {\frak m}_2 \ni t .
\]
Thus, $t^n$ is contained in $I$ for a sufficiently large $n$.
In particular, $t^n$ is in ${\frak m}/{\frak p}R$.

Consider the following commutative diagram:
\[
\begin{array}{ccc}
K[t^n]_{(t^n)} & \longrightarrow & R/{\frak p}R \\
\downarrow & & \downarrow \\
K[t]_{(t)} & \longrightarrow & \overline{S/{\frak p}} \otimes_{S/{\frak p}} R/{\frak p}R
\end{array}
\]
The morphism in the left-hand-side, as well as the bottom one
is finite. Hence, all morphisms are finite.

Put $s = t^n$.
Then, $s$ satisfies all the requirements.
\qed

Let $R \stackrel{\xi}{\longrightarrow} R/{\frak p}R$ be 
the natural surjective morphism.
We set $A = \xi^{-1} (K[s]_{(s)})$.
\[
\begin{array}{ccc}
R & \stackrel{\xi}{\longrightarrow} & R/{\frak p}R \\
\uparrow & \Box & \uparrow \\
A & \rightarrow & K[s]_{(s)}
\end{array}
\]
In the rest of this section, we shall prove that 
the ring $A$ satisfies the required condition.

Next we shall prove the following:

\begin{Claim}\label{claim2}
The morphism $A \rightarrow R$ is finite birational, and
$A$ is essentially of finite type over $K$ of dimension $d +1$. 
\end{Claim}

\proof
Remark that
\[
A \supset {\frak p}R \neq 0
\]
since the dimension of $R$ is at least $2$.
Take $0 \neq a \in {\frak p}R$. 
Since $A[a^{-1}] = R[a^{-1}]$,
$A \rightarrow R$ is birational.

One can prove that $A$ is a Noetherian ring by Eakin-Nagata's theorem.
However, here, we shall prove that $A$ is essentially of finite type
over $K$ without using Eakin-Nagata's theorem.

Let $B$ be the integral closure of $K[s]$ in $R/{\frak p}R$.
Remark that $B$ is of finite type over $K$.

Since $R/{\frak p}R$ is finite over $K[s]_{(s)}$,
$B \otimes_{K[s]} K[s]_{(s)} = R/{\frak p}R$.
\[
\begin{array}{ccccc}
R & \stackrel{\xi}{\longrightarrow} & R/{\frak p}R & \longleftarrow & B \\
\uparrow & & \uparrow & & \uparrow \\
S & & K[s]_{(s)} & \longleftarrow & K[s] 
\end{array}
\]
Take an element $s' \in R$ that satisfies $\xi(s') = s$.
Suppose $S = K[s_1, \ldots, s_n]$.
Since $B \otimes_{K[s]} K[s]_{(s)} = R/{\frak p}R$,
there exist $g_i \in B$ and $f_i \in K[s] \setminus (s)$ 
such that $\xi(s_i) = g_i/f_i$ for $i = 1, \ldots, n$.
Take an element $f'_i \in K[s']$ such that $\xi (f'_i) = f_i$
for $i = 1, \ldots, n$.
Put
\[
S' = K[s', s_1f'_1, \ldots, s_nf'_n] .
\]
Remark that $R$ is a localization of $S'$, and
$\xi(S') \subset B$.
Since $B$ is of finite type over $K$,
there exists a ring $D$ that satisfies
\begin{itemize}
\item
$S' \subset D \subset R$
\item
$D$ is of finite type over $K$,
\item
$R$ is a localization of $D$, and
\item
$\xi(D) = B$.
\end{itemize}
Put $\phi = \xi|_D$ and $E = \phi^{-1}(K[s])$.
Then, $D$ is finite over $E$.
\[
\begin{array}{ccc}
D & \stackrel{\phi}{\longrightarrow} & B \\
\uparrow & \Box & \uparrow \\
E & \rightarrow & K[s]
\end{array}
\]
Since $B \otimes_{K[s]} K[s]_{(s)} = R/{\frak p}R$,
there is only one prime ideal $N$ of $B$ lying over $(s) \subset K[s]$.
Therefore, $\phi^{-1}(N)$ is the only one prime ideal 
lying over the prime ideal $\phi^{-1}((s))$ of $E$.
Localizing all the rings in the above diagram, we have the following diagram:
\[
\begin{array}{ccc}
D \otimes_EE_{\phi^{-1}((s))} & \longrightarrow & 
B \otimes_EE_{\phi^{-1}((s))} \\
\uparrow & \Box & \uparrow \\
E_{\phi^{-1}((s))} & \rightarrow & K[s] \otimes_EE_{\phi^{-1}((s))}
\end{array}
\]
Remark that $D \otimes_EE_{\phi^{-1}((s))} = R$, 
$K[s] \otimes_EE_{\phi^{-1}((s))} = K[s]_{(s)}$ and
$B \otimes_EE_{\phi^{-1}((s))} = R/{\frak p}R$.
Therefore, $A$ coincides with $E_{\phi^{-1}((s))}$.

Since $D$ is finite over $E$ and $D$ is of finite type over $K$,
$E$ is also of finite type over $K$.

Therefore, we know that $A$ is essentially of finite type over $K$
and $R$ is finite over $A$.
It is easy to see
\[
\dim A = \dim R = \dim S = d+1 .
\]
\qed

In particular, 
$A$ is a homomorphic image of a regular local ring $T$.
Therefore, we have the commutative diagram
\[
\begin{array}{ccc}
G_0(A)_{\Bbb Q} & \stackrel{\tau_{A/T}}{\longrightarrow} & 
A_*(A)_{\Bbb Q} \\
\downarrow & & \downarrow \\
G_0(\widehat{A})_{\Bbb Q} & \stackrel{\tau_{\widehat{A}/\widehat{T}}}{\longrightarrow} & 
A_*(\widehat{A})_{\Bbb Q}
\end{array}
\]
by the singular Riemann-Roch theorem (Chapter~18, 20 in \cite{F}).
Remark that the horizontal maps in the above diagram are isomorphisms.
Therefore, in order to prove that 
${\rm Ker}(G_0(A)_{\Bbb Q} \rightarrow G_0(\widehat{A})_{\Bbb Q})$
is not $0$,
it is sufficient to prove that 
${\rm Ker}(A_1(A)_{\Bbb Q} \rightarrow A_1(\widehat{A})_{\Bbb Q})$ 
is not $0$.

The diagram
\[
\begin{array}{ccc}
R & \longrightarrow & \widehat{R} \\
\uparrow &  & \uparrow \\
A  & \longrightarrow & \widehat{A} 
\end{array}
\]
induces the commutative diagram
\begin{equation}\label{diag1}
\begin{array}{ccc}
A_1(R)_{\Bbb Q} & \longrightarrow & A_1(\widehat{R})_{\Bbb Q} \\
\downarrow &  & \downarrow \\
A_1(A)_{\Bbb Q} & \longrightarrow & A_1(\widehat{A})_{\Bbb Q} 
\end{array}
\end{equation}
where the vertical maps are induced by the finite morphisms
$A \rightarrow R$ and $\widehat{A} \rightarrow \widehat{R}$,
and the horizontal maps are induced by the completions
$A \rightarrow \widehat{A}$ and $R \rightarrow \widehat{R}$.

The top map in the diagram~(\ref{diag1}) is an isomorphism 
by assumption~1 of Lemma~\ref{lem1}.

Here we shall show, for each prime ideal of $A$, 
there exists only one prime ideal 
of $R$ lying over it.
Let ${\frak q}$ be a prime ideal of $A$.
Recall that the conductor ideal ${\frak p}R$ is a prime ideal
of both $A$ and $R$.
If ${\frak q}$ does not contain ${\frak p}R$, then $A_{\frak q}$ coincides with
$R \otimes_AA_{\frak q}$.
Therefore there exists only one prime ideal of $R$
lying over ${\frak q}$ in this case.
Next suppose that ${\frak q}$ contains ${\frak p}R$.
Then ${\frak q}$ is either ${\frak p}R$ or the unique maximal ideal of $A$.
In any cases, there exists only one prime ideal of $R$
lying over ${\frak q}$.

Consider the following commutative diagram:
\[
\begin{array}{ccccccccc}
0 & \longrightarrow & {\mbox {\it Rat}}_1(R) & \longrightarrow & Z_1(R) & \longrightarrow & A_1(R) & \longrightarrow & 0 \\
& & \downarrow & & \downarrow & & \downarrow & & \\
0 & \longrightarrow & {\mbox {\it Rat}}_1(A) & \longrightarrow & Z_1(A) & \longrightarrow & A_1(A) & \longrightarrow & 0
\end{array}
\]
We refer the reader to Chapter~1 in \cite{F} for definition
of ${\it Rat}_*$ and $Z_*$.
Since the morphism $A \rightarrow R$ is finite injective,
the cokernel of 
${\mbox {\it Rat}}_1(R) \rightarrow {\mbox {\it Rat}}_1(A)$ is 
torsion by Proposition~1.4 of Chapter~1 in \cite{F}.
Since, for each prime ideal of $A$, 
there is only one prime ideal of $R$ lying over it,
the map $Z_1(R) \rightarrow Z_1(A)$ is injective and 
the cokernel of it is a torsion module ${\Bbb Z}/(2v)$,
where $v = v_1(s) = v_2(s)$.
Therefore the map in the left-hand-side 
in diagram~(\ref{diag1}) is also an isomorphism.

By the commutativity of diagram~(\ref{diag1}),
we know that, in order to prove that 
${\rm Ker}(A_1(A)_{\Bbb Q} \rightarrow A_1(\widehat{A})_{\Bbb Q})$ 
is not $0$,
it is sufficient to show that
\[
{\rm Ker}(A_1(\widehat{R})_{\Bbb Q} 
\rightarrow A_1(\widehat{A})_{\Bbb Q}) = {\Bbb Q} .
\]

Since $\widehat{A}/({\frak p}R)\widehat{A} = \widehat{K[s]_{(s)}} = K[[s]]$,
$({\frak p}R)\widehat{A}$ is a prime ideal of $\widehat{A}$ 
of height $d$.
We have the following bijective correspondences:
\begin{eqnarray*}
& & \mbox{the set of prime ideals of $\widehat{R}$ 
lying over $({\frak p}R)\widehat{A}$ } \\
& \longleftrightarrow &
\mbox{the set of minimal prime ideals of $\widehat{R/{\frak p}R}$} \\
& \longleftrightarrow &
\mbox{the set of maximal ideals of $\overline{S/{\frak p}} \otimes_{S/{\frak p}} R/{\frak p}R$} \\
& \longleftrightarrow &
\{ P_1, P_2 \} ,
\end{eqnarray*}
where $\overline{S/{\frak p}} \otimes_{S/{\frak p}} R/{\frak p}R$ is the normalization
of $R/{\frak p}R$.
Therefore, there are just two prime ideals of $\widehat{R}$ 
lying over $({\frak p}R)\widehat{A}$.
We denote them by ${\frak p}_1$ and ${\frak p}_2$.

It is easy to see that ${\frak p}R$ is the conductor ideal of the ring 
extension
$A \subset R$, that is,
\[
{\frak p}R = A :_AR .
\]
Then, $({\frak p}R)\widehat{A} = \widehat{A} :_{\widehat{A}} \widehat{R}$ is
satisfied.
Therefore, $({\frak p}R)\widehat{A}$ is the conductor ideal of the ring 
extension
$\widehat{A} \subset \widehat{R}$.
Consider the map
\[
\varphi : Z_1(\widehat{R}) \longrightarrow Z_1(\widehat{A}) .
\]
Let ${\frak q}$ be a prime ideal of $\widehat{A}$ of height $d$.
If ${\frak q}$ does not contain the conductor ideal $({\frak p}R)\widehat{A}$, 
then there exists only one prime ideal ${\frak q}'$ of $\widehat{R}$
lying over ${\frak q}$. 
Furthermore, $\widehat{A}/{\frak q}$ is birational to $\widehat{R}/{\frak q}'$.
Therefore,
\[
\varphi ([\spec \widehat{R}/{\frak q}' ]) = [\spec \widehat{A}/{\frak q} ] .
\]
Here, 
we shall show
\[
\varphi ([\spec \widehat{R}/{\frak p}_1 ])
= \varphi ([\spec \widehat{R}/{\frak p}_2 ])
= v [\spec \widehat{A}/({\frak p}R)\widehat{A} ] ,
\]
where $v = v_1(s) = v_2(s)$.
Recall that
\[
\widehat{O_{v_1}} \times \widehat{O_{v_2}}
= (\overline{R/{\frak p}R})^\wedge
= \overline{\widehat{R}/{\frak p}\widehat{R}}
= \overline{\widehat{R}/{\frak p}_1}
\times \overline{\widehat{R}/{\frak p}_2} .
\]
Therefore, we may assume 
$\widehat{O_{v_i}} \simeq \overline{\widehat{R}/{\frak p}_i}$ 
for $i = 1, 2$.
Then, we have
\begin{eqnarray*}
& & \rank_{\widehat{A}/({\frak p}R)\widehat{A}} \widehat{R}/{\frak p}_i
= \rank_{\widehat{A}/({\frak p}R)\widehat{A}} \overline{\widehat{R}/{\frak p}_i}
= \rank_{\widehat{A}/({\frak p}R)\widehat{A}} \widehat{O_{v_i}} 
= \rank_{K[[s]]} \widehat{O_{v_i}} \\
& = & \dim_K \widehat{O_{v_i}}/s\widehat{O_{v_i}}
= \dim_K O_{v_i}/sO_{v_i} = v
\end{eqnarray*}
for $i = 1, 2$.
Here, $\dim_K$ means the dimension of the given $K$-vector space.

Thus, we have the following exact sequence
\[
0 \longrightarrow {\Bbb Z} \cdot 
([\spec \widehat{R}/{\frak p}_1 ] - [\spec \widehat{R}/{\frak p}_2 ])
\longrightarrow 
Z_1(\widehat{R})
\longrightarrow 
Z_1(\widehat{A})
\longrightarrow 
{\Bbb Z}/(v)
\longrightarrow
0 .
\]
Consider the following diagram:
\[
\begin{array}{ccccccccc}
0 & \longrightarrow & {\mbox {\it Rat}}_1(\widehat{R}) & \longrightarrow & 
Z_1(\widehat{R}) & \longrightarrow & A_1(\widehat{R}) & \longrightarrow & 0 \\
& & \downarrow & & \downarrow & & \downarrow & & \\
0 & \longrightarrow & {\mbox {\it Rat}}_1(\widehat{A}) & \longrightarrow & 
Z_1(\widehat{A}) & \longrightarrow & A_1(\widehat{A}) & \longrightarrow & 0
\end{array}
\]
Since the morphism $\widehat{A} \rightarrow \widehat{R}$ is
finite injective,
the cokernel of  ${\mbox {\it Rat}}_1(\widehat{R}) \rightarrow {\mbox {\it Rat}}_1(\widehat{A})$ is torsion (c.f.\ Proposition~1.4 in \cite{F}).
Thus, we have the following exact sequence
\[
0 \longrightarrow {\Bbb Q} \cdot 
([\spec \widehat{R}/{\frak p}_1 ] - [\spec \widehat{R}/{\frak p}_2 ])
\longrightarrow 
A_1(\widehat{R})\subq
\longrightarrow 
A_1(\widehat{A})\subq
\longrightarrow 
0 .
\]
Therefore, we have only to prove
\[
[\spec \widehat{R}/{\frak p}_1 ] - [\spec \widehat{R}/{\frak p}_2 ]
\neq 0
\]
in $A_1(\widehat{R})\subq$.

Let $\widehat{\pi} : 
\widehat{Y} \rightarrow \spec \widehat{R} $ be the
blow-up at ${\frak m} \widehat{R}$.
Since $\widehat{\pi}^{-1}({\frak m} \widehat{R}) \simeq X$,
\[
A_1(X)_{\Bbb Q} \stackrel{i_*}{\rightarrow} A_1(\widehat{Y})_{\Bbb Q}
\stackrel{\widehat{\pi}_*}{\rightarrow} A_1(\widehat{R})_{\Bbb Q} \rightarrow 0
\]
is exact and
\[
\widehat{\pi}_* \left(
[\spec \overline{\widehat{R}/{\frak p}_1} ] - 
[\spec \overline{\widehat{R}/{\frak p}_2} ]
\right)
= 
[\spec \widehat{R}/{\frak p}_1 ] - 
[\spec \widehat{R}/{\frak p}_2 ] ,
\]
where $i : X \rightarrow \widehat{Y}$ is the inclusion.
Consider the following commutative diagram:
\[
\begin{array}{ccccc}
P_i & \longrightarrow & \{ P_1, P_2 \} & \longrightarrow & X \\
\downarrow & \Box & \downarrow & \Box & \downarrow \\
\spec O_{v_i}  & \longrightarrow & \spec \overline{R/{\frak p}} 
& \longrightarrow &  Y \\
& & \downarrow & & \downarrow \\
& & \spec R/{\frak p}  & \longrightarrow & \spec R 
\end{array}
\]
Take the fibre product with $\spec \widehat{R} $ over $\spec R $.
We may assume that $\spec \overline{\widehat{R}/{\frak p}_i} $
coincides with $\spec \widehat{O_{v_i}} $
for $i = 1, 2$ so that the following diagram commutes:
\[
\begin{array}{ccccccccc}
& & P_i & = & P_i & \longrightarrow & 
\{ P_1, P_2 \} & \longrightarrow & X \\
& & \downarrow & \Box & \downarrow & \Box & \downarrow & \Box 
& \downarrow \\
\spec \overline{\widehat{R}/{\frak p}_i}  & = & 
\spec \widehat{O_{v_i}}  & \longrightarrow & 
\spec O_{v_i} \otimes_R \widehat{R}  &
\longrightarrow & \spec \overline{R/{\frak p}}\otimes_R \widehat{R} 
& \longrightarrow &  \widehat{Y} \\
\downarrow & & & & & & \downarrow & & \downarrow \\
\spec \widehat{R}/{\frak p}_i  & 
\multicolumn{5}{c}{\longrightarrow} & 
\spec \widehat{R}/{\frak p}\widehat{R} 
& \longrightarrow & \spec \widehat{R} 
\end{array}
\]
Assume that
\[
[\spec \widehat{R}/{\frak p}_1 ] - [\spec \widehat{R}/{\frak p}_2 ] = 0
\]
in $A_1(\widehat{R})\subq$.
Then, there exists $\delta \in A_1(X)_{\Bbb Q}$ such that
\[
i_*(\delta) = 
[\spec \overline{\widehat{R}/{\frak p}_1} ] - 
[\spec \overline{\widehat{R}/{\frak p}_2} ] .
\]
Here, consider the map
\[
A_1(\widehat{Y})_{\Bbb Q} \stackrel{i^!}{\rightarrow} 
A_0(X)_{\Bbb Q} , 
\]
that is taking the intersection with 
$\widehat{\pi}^{-1}({\frak m} \widehat{R}) = X$.
Since $i^! i_*(\delta) = -h \cdot \delta$ and
\[
i^! \left( 
[\spec \overline{\widehat{R}/{\frak p}_1} ] - 
[\spec \overline{\widehat{R}/{\frak p}_2} ] \right) 
=
i^! \left( 
[\spec \widehat{O_{v_1}} ] - 
[\spec \widehat{O_{v_2}} ] \right) 
= 
[P_1] - [P_2] ,
\]
we have 
\[
[P_1] - [P_2] = -h \cdot \delta .
\]
It contradicts to 
\[
[P_1] - [P_2] \neq 0
\]
in $A_0(X)_{\Bbb Q}/h \cdot A_1(X)_{\Bbb Q}$.

We have completed the proof of Lemma~\ref{lem1}.

\section{A proof of Lemma~\ref{lem2}}
We shall give a proof of Lemma~\ref{lem2}
in this section.

Suppose that $S = {\Bbb C}[x_0,x_1,x_2]/(f)$ and
$X = \proj S $ satisfy the assumption in Lemma~\ref{lem2}.
Let $Z$ be the projective cone of $X$, that is,
$Z =  \proj {\Bbb C}[x_0,x_1,x_2,x_3]/(f) $.

Let $W \stackrel{\xi}{\rightarrow} Z$ be the blow-up at $(0,0,0,1)$.
We set $X_\infty = V_+(x_3)$ and $X_0 = \xi^{-1}((0,0,0,1))$.
Remark that both of $X_0$ and $X_\infty$ are isomorphic to $X$.
Then, $W \stackrel{\eta}{\rightarrow} X$ is a ${\Bbb P}^1$-bundle.

Take any two closed points $Q_1, Q_2 \in X$.
We set $L_i = \eta^{-1}(Q_i)$ for $i = 1, 2$.
Consider the Weil divisor
$L_1 + L_2 + X_\infty$ on $W$.
Here we shall prove the following:

\begin{Claim}\label{claim3}
The complete linear system 
$| L_1 + L_2 + X_\infty |$ is base-point free, 
and the induced morphism
$W \stackrel{f}{\rightarrow} {\Bbb P}^n$ satisfies that
$\dim f(W) \geq 2$.
\end{Claim}

\proof
Since the complete linear system $| Q_1 + Q_2 |$ on $X$ is
base-point free,
so is $| L_1 + L_2 |$.
Since the complete linear system $| X_\infty |$ is
base-point free,
so is $| L_1 + L_2 + X_\infty |$.

In order to show $\dim f(W) \geq 2$,
we have only to show that the set 
\[
\left\{ a \in R(W)^\times \mid 
{\rm div}(a) + L_1 + L_2 + X_\infty \geq 0 \right\}
\]
contains two algebraically independent elements over ${\Bbb C}$.

Note that, since $W \stackrel{\eta}{\rightarrow} X$ 
is a surjective morphism, $R(X)$ is contained in $R(W)$.
Consider
\begin{eqnarray*}
H^0(W, {\cal O}_W(L_1 + L_2 + X_\infty))
& = & 
\left\{ a \in R(W)^\times \mid 
{\rm div}(a) + L_1 + L_2 + X_\infty \geq 0 \right\} 
\cup \{ 0 \} \\
H^0(X, {\cal O}_X(Q_1 + Q_2)) 
& = & 
\left\{ a \in R(X)^\times \mid 
{\rm div}(a) + Q_1 + Q_2 \geq 0 \right\} 
\cup \{ 0 \} .
\end{eqnarray*}
It is easy to see 
\[
H^0(W, {\cal O}_W(L_1 + L_2 + X_\infty))
\supset
H^0(X, {\cal O}_X(Q_1 + Q_2)) \supset {\Bbb C} .
\]
The set $H^0(X, {\cal O}_X(Q_1 + Q_2))$
contains a transcendental element over ${\Bbb C}$.
Since $R(X)$ is algebraically closed in $R(W)$ and
\[
H^0(W, {\cal O}_W(L_1 + L_2 + X_\infty))
\neq
H^0(X, {\cal O}_X(Q_1 + Q_2)) ,
\]
$H^0(W, {\cal O}_W(L_1 + L_2 + X_\infty))$
contains two algebraically independent elements over ${\Bbb C}$.
\qed

Since $| L_1 + L_2 + X_\infty |$ is base-point free
as in Claim~\ref{claim3},
\[
{\rm div}(a) + L_1 + L_2 + X_\infty
\]
is smooth for a general element $a \in 
H^0(W, {\cal O}_W(L_1 + L_2 + X_\infty)) \setminus \{ 0 \}$
(e.g., \RMN{3} Corollary~10.9 in \cite{Ha}).
Since $\dim f(W) \geq 2$
as in Claim~\ref{claim3},
\[
{\rm div}(a) + L_1 + L_2 + X_\infty
\]
is connected for any $a \in 
H^0(W, {\cal O}_W(L_1 + L_2 + X_\infty)) \setminus \{ 0 \}$
(e.g., \RMN{3} Exercise~11.3 in \cite{Ha}).

Let $\{ a_1, \cdots, a_n \}$ be a ${\Bbb C}$-basis
of $H^0(W, {\cal O}_W(L_1 + L_2 + X_\infty))$.
Let $\alpha_i$ be the local equation defining
the Cartier divisor ${\rm div}(a_i) + L_1 + L_2 + X_\infty$
for $i = 1, \ldots, n$.
For $c = (c_1, \ldots, c_n) \in {\Bbb C}^n \setminus \{ (0, \ldots, 0) \}$,
$D_c$ denotes the Cartier divisor on $W$ defined by 
$c_1\alpha_1 + \cdots + c_n \alpha_n$.

For a general point $c \in {\Bbb C}^n$,
$D_c$ does not contain $X_0$ as a component and
$D_c$ intersect with $X_0$ at two distinct points.
Recall that $X_0$ is isomorphic to $X$.
Set $D_c \cap X_0 = \{ Q_{c1}, Q_{c2} \} \subset X$.

Choose $e \in X$ such that the Weil divisor $3e$ coincides with
the very ample divisor corresponding to 
the embedding $X = \proj S $.
We regard the set of closed points of the elliptic curve $X$ 
as an abelian group with unit $e$ as in the usual way.

Let $\varphi : X \rightarrow {\Bbb P}^1_{\Bbb C}$ be the 
morphism defined by $| 2 e |$.

For a general point $c \in {\Bbb C}^n$,
we set
\[
\theta(c) = \varphi (Q_{c1} \ominus Q_{c2}) 
\in {\Bbb P}^1_{\Bbb C} ,
\]
where $\ominus$ means the difference in the group $X$.
One can prove that there exists a non-empty Zariski open set
$U$ of ${\Bbb C}^n$ such that 
$\theta|_U : U \rightarrow {\Bbb P}^1_{\Bbb C}$ is 
a non-constant morphism and $D_c$ is smooth 
connected for any $c \in U$.
Then, there exists a a non-empty Zariski open set of ${\Bbb P}^1_{\Bbb C}$
which is contained in ${\rm Im}(\theta|_U)$.
Let $F$ be the set of elements of $X$ of order finite.
Then, it is well-known that $F$ 
is a countable set.
In particular, $\varphi(F)$ does not contain ${\rm Im}(\theta|_U)$.
Therefore, there exists $c \in U$
such that $\theta(c) \not\in \varphi(F)$.
Then, $D_c$ is a smooth connected curve in $W$
such that $D_c$ intersect with $X_0 \simeq X$
at two points $\{ P_1, P_2 \}$ transversally
such that $P_1 \ominus P_2$ has order infinite
in $X$.

Let $\phi : X \rightarrow A_0(X)$ be a map
defined by $\phi(P) = [P] - [e]$.
It is well known that $\phi$ is a group homomorphism.
We have the following exact sequence:
\[
0 \longrightarrow X \stackrel{\phi}{\longrightarrow}
A_0(X) \stackrel{{\rm deg}}{\longrightarrow}
{\Bbb Z} \longrightarrow 0 
\]
Since ${\rm deg}(h) = 3$, we have an isomorphism
\[
X\otimes_{\Bbb Z}{\Bbb Q} \stackrel{\overline{\phi}}{\simeq} A_0(X)\subq/h A_1(X)\subq .
\]
By definition, we have 
\[
0 \neq \overline{\phi} (P_1\ominus P_2)
= [P_1] - [P_2]
\]
in $A_0(X)\subq/h A_1(X)\subq$.

Let $Y$ be the blow-up of $\spec S $ at the origin.
Then, $Y$ is an open subvariety of $W$.
We set $C = D_c \cap Y$.
Then $C$ satisfies assumption~2 in Lemma~\ref{lem2}.

Since $H^1(X, {\cal O}_X(n)) = 0$ for $n > 0$,
we have ${\rm Cl}(R) \simeq {\rm Cl}(\widehat{R})$ 
by Danilov's Theorem (Corollary in 497p and Proposition~8
in \cite{D}).
Therefore, $R$ satisfies 
assumption~1 in Lemma~\ref{lem1}.

We have completed the proof the Lemma~\ref{lem2}.

\begin{Remark}
\begin{rm}
Let $A$ be a $2$-dimensional local ring constructed using
Lemma~\ref{lem1} and Lemma~\ref{lem2}.
Since $A$ and $\widehat{A}$ are $2$-dimensional excellent local domains,
we have the following isomorphisms:
\begin{eqnarray*}
G_0(A) & \simeq & {\Bbb Z} \oplus A_1(A) \\
G_0(\widehat{A}) & \simeq & {\Bbb Z} \oplus A_1(\widehat{A}) 
\end{eqnarray*}
Therefore, 
\[
{\rm Ker}(G_0(A) \rightarrow G_0(\widehat{A})) 
\simeq {\rm Ker}(A_1(A) \rightarrow A_1(\widehat{A}) ) .
\]
Using it, we can prove that
\[
{\rm Ker}(G_0(A) \rightarrow G_0(\widehat{A})) 
\simeq {\Bbb Z} 
\]
as follows.
Consider the following diagram
\[
\begin{array}{ccc}
 & & 0 \\
 & & \downarrow \\
0 & & {\Bbb Z} \\
\downarrow & & \downarrow \\
A_1(R) & \specialarrow{\longrightarrow}{f}{\displaystyle \sim} & 
A_1(\widehat{R}) \\
\phantom{i} \downarrow i & & \downarrow \\
A_1(A) & \stackrel{g}{\longrightarrow} & A_1(\widehat{A}) \\
\downarrow & & \downarrow \\
{\Bbb Z}/(2v) & \longrightarrow & {\Bbb Z}/(v) \\
\downarrow & & \downarrow \\
0 & & 0 
\end{array}
\]
Let $\alpha_i$ be the element of $A_1(R)$
such that $f(\alpha_i) = [\spec \widehat{R}/{\frak p}_i]$
for $i = 1, 2$.
Then, the kernel of $g$ is generated by
\[
i(\alpha_1) - v [\spec A/{\frak p}R] .
\]
Here, note that
\[
2 \left( i(\alpha_1) - v [\spec A/{\frak p}R] \right)
= i(\alpha_1) - i(\alpha_2) .
\]
Since the kernel of $g$ is not torsion,
it must be isomorphic to ${\Bbb Z}$.
\end{rm}
\end{Remark}

\noindent
\begin{tabular}{l}
Kazuhiko Kurano \\
Department of Mathematics \\
Faculty of Science and Technology \\
Meiji University \\
Higashimita 1-1-1, Tama-ku \\
Kawasaki 214-8571, Japan \\
{\tt kurano@math.meiji.ac.jp} \\
{\tt http://www.math.meiji.ac.jp/\~{}kurano}
\end{tabular}

\vspace{2mm}

\noindent
\begin{tabular}{l}
Vasudevan Srinivas \\
School of Mathematics \\
Tata Institute of Fundamental Research \\
Homi Bhabha Road, Mumbai-400005 \\
India \\
{\tt srinivas@math.tifr.res.in} 
\end{tabular}


\begin{thebibliography}{99}

\bibitem{D}  {\sc V. I. Danilov},
{\it The group of ideal classes of a complete ring}, 
Math.\ USSR Sbornik {\bf 6} (1968), 493--500.

\bibitem{Dao} {\sc H. Dao},
{On injectivity of maps between Grothendieck groups induced by completion},
preprint

\bibitem{F} {\sc W.~Fulton},
{\it Intersection Theory, 2nd Edition},
Springer-Verlag, Berlin, New York, 1997.

\bibitem{GS} {\sc H. Gillet and C. Soul\'{e}},
{\it K-th\'{e}orie et nullit\'{e} des multiplicites d'intersection},
C. R. Acad.\ Sci.\ Paris Ser.\ I Math. {\bf 300} (1985), 71--74.

\bibitem{Ha} {\sc R. Hartshorne},
{\it Algebraic Geometry,} Graduate Texts in Math., No.\ 52,
Springer-Verlag, Berlin and New York, 1977.

\bibitem{H} {\sc M. Hochster},
{\it Thirteen open questions in Commutative Algebra},
talk given at LipmanFest,
July 2004, available online at 
{\tt http://www.math.lsa.umich.edu/hochster/Lip.test.pdf}




\bibitem{KK} {\sc Y. Kamoi and K. Kurano},
{\it On maps of Grothendieck groups induced by completion},
J.\ Algebra {\bf 254} (2002), 21--43.


\bibitem{K11} {\sc K. Kurano},
{\it A remark on the Riemann-Roch formula for affine schemes
associated with Noetherian local rings}, 
T\^ohoku Math.\ J. {\bf 48} (1996), 121--138.


\bibitem{K16} {\sc K. Kurano},
{\it On Roberts rings}, J. Math. Soc. Japan {\bf 53} (2001), 333--355.


\bibitem{K23} {\sc K. Kurano},
{\it Numerical equivalence defined 
on Chow groups of Noetherian local rings},
Invent. Math., {\bf 157} (2004), 575--619.






\bibitem{R2} {\sc P. C. Roberts},
{\it The vanishing of intersection multiplicities and perfect complexes},
Bull.\ Amer.\ Math.\ Soc. {\bf 13} (1985), 127--130.







\end{thebibliography}
\end{document}